\documentclass{aptpub}
\usepackage{amsmath}
\usepackage{latexsym}
\usepackage{amsmath}
\usepackage{amssymb}
\usepackage{amsfonts}
\usepackage{epsfig}
\usepackage{times}

\authornames{Lu, Radovanovi\'c} 

\shorttitle{Asymptotic Blocking Probabilities in the Loss Networks
with Subexponential Demands}

\def\qed{\hfill{$\Diamond$} \\}

\def\expect{{\mathbb  E}}
\def\ex{{\mathbb  E}}

\def\Pr{{\mathbb P}}
\def\real{{\mathbb  R}}
\def\nat{{\mathbb  N}}

\def\proof{\noindent{\bf Proof:} }
\def\eqdef{\triangleq}
\def\eqd{\stackrel{d}{=}}
\def\leqd{\stackrel{d}{\le}}

\def\th {\mit\theta}

\renewcommand{\baselinestretch}{2}

\begin{document}

\title{Asymptotic Blocking Probabilities in Loss Networks with\\
Subexponential Demands}

\vspace{0.3in}
\authorone[IBM T.J.\ Watson Research Center]{Yingdong Lu and Ana Radovanovi\'c}
\addressone{Mathematical Sciences Department,
IBM Thomas J.\ Watson Research Center, Yorktown Heights, NY 10598,
\{yingdong, aradovan\}@us.ibm.com.}

\begin{abstract}

The analysis of stochastic loss networks has long been of interest
in computer and communications networks and is becoming important in
the areas of service and information systems. In traditional
settings, computing the well known Erlang formula for blocking
probability in these systems becomes intractable for larger resource
capacities. Using compound point processes to capture stochastic
variability in the request process, we generalize existing models in
this framework and derive simple asymptotic expressions for blocking
probabilities. In addition, we extend our model to incorporate
reserving resources in advance. Although asymptotic, our experiments
show an excellent match between derived formulas and simulation
results even for relatively small resource capacities and relatively
large values of blocking probabilities.

\vspace{10mm} \noindent {\bf Keywords:} loss networks;
subexponential distributions.
\end{abstract}

\newpage

\section{Introduction}

The problem of satisfying a stream of customer (user) requirements
from resources of finite capacities for some random processing
time has long been present in many areas such as telephone and
communication networks, inventory control (rental industry) and,
recently, workforce management. For all of these applications,
system dynamics can be described as follows. Requests for
resources arrive according to some point process in time. If there
are enough available (non-engaged) resources to satisfy their
requirements at the moment of arrival, required resources are
committed for some random time that represents their processing
duration (holding time) after which they are released and become
available to accommodate future requests. In the case of
insufficient amount of available resources at the moment of its
arrival, a request is lost. The previously described system is
usually referred to as a {\em loss network}, and one of the
commonly analyzed performance metrics is the blocking probability,
i.e., probability that an incoming request is lost due to
insufficient amount of available resources to satisfy its
requirements.

Loss networks with fixed resource requirements have been
intensively analyzed in the context of circuit-switched networks.
Let requests require resources of $K<\infty$ different types for
some random generally distributed processing time with finite
mean. Furthermore, assume that requests belong to $M$ different
classes characterized by their resource requirements, processing
durations, arrival rates. Then, assuming that requests of
different types arrive according to mutually independent Poisson
processes, by PASTA property (\cite{WOL89}), blocking probability
$B_l$ of an incoming request of type $1\le l\le M$ is equal to the
sum of probabilities of blocking states for $l$ type request and
is computed using the generalized Erlang formula (e.g., see
\cite{KEL91s}), i.e.,
\begin{eqnarray*}
B_l=1-G({\bf C})^{-1} G({\bf C} -{\bf A} {\bf e}_l),
\end{eqnarray*}
where
\begin{eqnarray*} G({\bf C})
=\left(\sum_{{\bf n} \in {\cal S}({\bf C})}\prod_{l=1}^M
\frac{\rho_l^{n_l}}{n_l!} \right)
\end{eqnarray*}
and
\begin{eqnarray}
\label{eqn:domain} {\cal S}({\bf C}):= \{ {\bf n}\in
\mathbb{Z}^M_+: {\bf A} {\bf n} \le {\bf C}\},
\end{eqnarray}
where ${\bf n}=(n_1,\dots ,n_{M})$ and ${\bf C}=(C_1,\dots,C_K)$.
In the previous expressions $C_k$, $1\le k\le K$, is capacity of
resource type $k$, ${\bf A}=[A_{kl}]$ is a $K\times M$ matrix,
where $A_{kl}$ represents the amount of resources of type $1\le
k\le K$ required by a request of type $1\le l\le M$, and $\rho_l$,
$1\le l\le M$, represent traffic intensities of $l$ type requests
(computed as $\rho_l=\lambda_l/\mu_l$, where $\lambda_l$ is the
arrival rate of $l$ type requests and $1/\mu_l$ is the
corresponding mean processing time). Furthermore, $e_l$ is a $M$
dimensional vector with the $l$th component equal to one and the
rest equal to zero. In the case of a single resource type and a
single request class with exponentially distributed processing
times, blocking probability was first expressed by Erlang in 1917
(see \cite{ERL17}). Later on, it was shown that the Erlang formula
holds under more general assumptions on call holding time
distributions (see \cite{SEV57}) and in the case of Poisson
arrivals with retrials (see \cite{BON06}). It is noteworthy to
point out the difference between the Erlang loss network and a
queue with finite buffer. The two systems follow very different
dynamics resulting in a different behavior and, therefore, their
analysis (e.g., see \cite{JEL97c} and \cite{ASP05}).

It is easy to see that the cardinality of the state space ${\cal
S}({\bf C})$ in (\ref{eqn:domain}) increases exponentially in the
norm of vector ${\bf C}$, i.e., $|{\bf C}|\equiv
\sum_{i=1}^{K}|C_i|$. It is shown in \cite{LOU94} that the
calculation of $G({\bf C})$ is a $\sharp P$-complete problem,
which belongs to a class of problems that are at least as hard as
$NP$-complete problems. To this end, many approximation techniques
for evaluating blocking probabilities in large loss networks have
been proposed. One of the most popular ones is known as Erlang
fixed point method. The main idea of this approximation is to
assume that deficiencies of different resource types happen
independently. The application of the Erlang fixed point method
can be traced back as early as 50's (e.g., see \cite{WIL56}). In
\cite{KEL86}, Kelly studied the performance of the Erlang fixed
point method and established its relation to a nonlinear
optimization problem. He also proved uniqueness of the fixed point
and its asymptotic exactness when resource capacities and arrival
rates grow with the same rate (see \cite{KEL91s}). Some of the
related practical aspects of Kelly's analysis were investigated in
\cite{WIT85}. The Erlang fixed point method is further refined in
\cite{ZAC91}. There are also many other types of approximations
such as recursive algorithm in \cite{KAUF81}, or unified approach
based on large deviations for all (light, critical and heavy)
traffic regimes in \cite{GLM93}. Overall, except from the bounds
in \cite{GLM93}, these methods make use of the structural
properties of the Erlang formula and, hence, largely rely on the
Poisson assumption for call arrivals. Another restriction of the
above models is that the amount of resource requirements are
assumed to be fixed; in fact, it is assumed that they are $(0,1)$
parameters in most of the cases considered. Meanwhile, we see in
many applications that resource requirements could be highly
variable and their distributions possibly long-tailed; for
specific examples, see \cite{HEL96}, \cite{JLS95a} and
\cite{KRM97}. Furthermore, more recently, loss networks models
have been applied in the context of workforce management
applications (see \cite{LU06}), where requests behavior is even
more volatile and extreme.

In this paper, we analyze loss networks that have renewal arrivals
and random resource requirements. In particular, we assume that
request arrivals follow a compound renewal process, with the
corresponding holding times being arbitrarily distributed with
finite mean, independent of each other and arrival points. In
order to cope with variability in resource requirements, we model
them as subexponential random variables. We obtain a simple and
explicit asymptotic expressions for blocking probabilities when
capacities of resources grow. For the case of a single resource
loss network, we show that the stationary blocking probability is
approximately equal to the tail of the resource requirement
distribution. In addition, we extend our results to allow advance
reservations of resources. Finally, we investigate general
(multiple resources and arbitrary topology) loss networks and show
that the asymptotic blocking probability behaves as the tail of
the heaviest-tailed resource requirement. Although asymptotic, our
numerical experiments show an excellent accuracy of the derived
formulas even for relatively small capacities and relatively large
values of blocking probabilities, suggesting wide applicability of
the obtained results.

Our paper is organized as follows. In Section \ref{sec:model}, we
introduce our model in the context of a single resource type. Then,
in Subsection \ref{sec:singlepool}, we state and prove our main
result in Theorem \ref{theorem:th1}, while in Subsection
\ref{sec:advancereservation}, we extend it to the case of advance
reservations. Further extension to the analysis of the stationary
blocking probability in the case of general loss networks is stated
and proved in Theorem \ref{theorem:th2} of Section
\ref{sec:network}. Our simulation experiments for some specific
cases of arrival processes and resource requirements are presented
in Section \ref{sec:experiment}. Finally, we conclude our paper in
Section \ref{sec:concl}. A discussion and the proof of existence of
the stationary blocking probability is presented in the Appendix.

\section{Systems with one resource type}
\label{sec:model}
Let requests for resources from a common resource pool of capacity
$C< \infty$ arrive at time points $\{\tau_{n},-\infty < n<
\infty\}$ that represent a renewal process with rate
$0<\lambda<\infty$, i.e., $\expect [\tau_n
-\tau_{n-1}]=1/\lambda$. At each point $\tau_{n}$, $B_{n}$ amount
of resources is requested. If available capacity is less than
$B_n$, this request is rejected (blocked); otherwise, it is
accepted and $B_n$ amount of resources will be occupied for the
length of time $\theta_{n}$. Sequences $\{B_{n}\}$ and
$\{\theta_n\}$ of i.i.d. random variables (r.v.) are assumed to be
mutually independent and independent of the arrival points
$\{\tau_{n}\}$; furthermore $\expect \theta_n <\infty$ for all
$n$. Let $B$ and $\theta$ denote random variables that represent
$\{B_{n}\}$, $\{\theta_n\}$, i.e., $\Pr[B>x]=\Pr[B_n>x],
\Pr[\th>y]=\Pr[\th_n>y]$, for any $n\in\mathbb{Z}$, $x\ge 0$ and
$y\ge 0$.

In this paper, we assume that $B$ is a subexponential random
variable, defined as follows (e.g., see \cite{GOK97}):
\begin{definition}Let $\{X_i\}$ be a sequence of positive i.i.d. random variables
with distribution function $F$ such that $F(x) <1$ for all $x>0$.
Denote by $\Bar {F}(x)=1-F(x)$, $x\ge 1$, the tail of $F$ and by
$\Bar {F}^{n\ast} = 1-F^{n\ast}(x)=\Pr [X_1+\dots +X_n
>x]$ the tail of the n-fold convolution of $F$. $F$ is
subexponential distribution function, denoted as  $F\in {\cal S}$,
if one of the following equivalent conditions holds:
\begin{itemize}
\item $\lim_{x\rightarrow \infty} \frac{\Bar {F}^{n\ast}(x)}{\Bar
{F} (x)}=n$ for some (all) $n\ge 2$, \item $\lim_{x\rightarrow
\infty} \frac{\Pr [X_1+\dots +X_n >x]}{\Pr [\max (X_1,\dots
,X_n)>x]}=1$ for some (all) $n\ge 2$.
\end{itemize}
\end{definition}
\vspace{0.1in}

For a brief introduction to subexponential distributions the
reader is referred to a recent survey \cite{GOK97}. This class of
distributions is fairly large and well known examples include
regularly varying (in particular Pareto), some Weibull, log-normal
and "almost" exponential distributions.

Next, let ${\cal N}^{(C)}_n$ be the set of indices $i< n$ of
resource requirements that arrive prior to $\tau_n$, are accepted,
and are {\em still active} by time $\tau_n$. Furthermore, let
$N^{(C)}_n\eqdef |{\cal N}^{(C)}_n|$ be a cardinality of set
${\cal N}^{(C)}_n$. Thus, the total amount of resources
$Q_n^{(C)}$ that an arrival at time $\tau_n$ finds engaged can be
expressed as $Q^{(C)}_{n}= \sum_{i\in {\cal N}^{(C)}_n} B_{i}$.

Our goal in this paper is to estimate the stationary blocking
probability, i.e.,
\begin{equation}
\label{eq:deflossrate} \Pr [Q^{(C)}_n +B_n>C],
\end{equation}
for large $C$. It can be shown that for the model introduced above
there exists a unique stationary distribution for $Q^{(C)}_n$ and,
therefore, the quantity in (\ref{eq:deflossrate}) is well defined.
The proof of this result is based on constructing a Markov chain
with general state space, of which $Q^{(C)}_n$ is a functional.
Then, by using a discrete version of Theorem 1 from \cite{SEV57},
we show that there exists a unique stationary distribution for the
constructed Markov chain (and, therefore, $Q^{(C)}_n$) which is
ergodic. Since this proof is not the main focus of this paper, we
present it in the Appendix.

In this paper we use the following standard notation. For any two
real functions $a(t)$ and $b(t)$ and fixed $t_{0}\in
\real\cup\{\infty\}$, let  $a(t)\thicksim b(t)$ as $t\rightarrow
t_{0}$ denote $\lim_{t\rightarrow t_{0}} [a(t)/b(t)]=1$.

\subsection{Blocking probability in a system with one resource type}
\label{sec:singlepool}

In this section we estimate the stationary blocking probability
$\Pr [Q^{(C)}_n+B_n>C]$ in a loss network with a single resource
pool when its capacity $C$ grows large.

\begin{theorem}
\label{theorem:th1} Let $\{B_n,-\infty<n<\infty\}$ be a sequence of
subexponential random variables with finite mean. Then, the
stationary blocking probability satisfies
\begin{equation}
\label{eq:asymmain} \Pr [Q^{(C)}_{n}+B_{n}>C]\sim \Pr
[B>C]\;\;\text{as}\;\;\text{$C\rightarrow \infty$.}
\end{equation}
\end{theorem}
\proof First, observe that a request will be lost if it requires
more than the total capacity $C$ and, therefore,
\begin{equation} \label{eq:lowbd5} \Pr [Q^{(C)}_n
+B_{n}>C]\ge \Pr [B>C]\;\;\text{ for all $C>0$.}
\end{equation}

In order to prove the asymptotic upper bound for $\Pr
[Q^{(C)}_n+B_{n}>C]$, we start by conditioning on the size of $B_n$
as
\begin{align}
\Pr [Q^{(C)}_{n}+B_{n}>C]&=\Pr [Q^{(C)}_{n}+B_{n}>C,
B_{n}>C]+\Pr [Q^{(C)}_{n}+B_{n}>C,B_{n}\le C]\nonumber\\
&\eqdef I_{1}+I_2.\label{eq:split1}
\end{align}
Note that $I_1$ is upper bounded by $\Pr [B>C]$. Next, we prove that
$I_2=o(\Pr [B>C])$ as $C\rightarrow \infty$. In view of the
definition of ${\cal N}^{(C)}_{n}$ from above,
\begin{equation}
\label{eq:crucial} I_2 = \Pr \left [\sum_{i\in {\cal N}^{(C)}_{n}}
B_i+B_{n}>C, B_{n}\le C \right].
\end{equation}
Observe that for $i\in {\cal N}^{(C)}_{n}$, $B_i$s are mutually
dependent which makes direct analysis of the expression in
(\ref{eq:crucial}) complex. For that reason, we sample the original
process of arrivals at points $\tau_i$ at which the requested amount
of resources $B_i$ is smaller or equal to $C$ and observe another
system of unlimited capacity with the sampled arrivals. Let ${\cal
N}_{n,s}$ be a set of request indices $i<n$ that belong to the
sampled process and are still active at time $\tau_n$, i.e.,
\begin{equation*}
{\cal N}_{s,n}=\{i<n|B_i\le C,\theta_i>\tau_n -\tau_i\}.
\end{equation*}
Note that the sampled process is renewal as well with rate
$\lambda \Pr [B\le C]/\Pr [B > C]$ and that resource requirements
$B_i$, $i\in {\cal N}_{s,n}$, are mutually independent.
Furthermore, since ${\cal N}^{(C)}_{n}\subset {\cal N}_{s,n}$, we
can upper bound $I_2$ in (\ref{eq:crucial}) by the probability
that the total amount of required resources in a new system
exceeds capacity $C$, i.e.,
\begin{equation}
\label{eq:refth11} I_2\le \Pr \left [\sum_{i\in {\cal N}_{s,n}}
B_i+B_{n}>C, B_{n}\le C \right].
\end{equation}
Now, in view of the results derived in \cite{EMG80} for every
integer $n$ and i.i.d. subexponential random variables $B_1,\dots
,B_n$, $\Pr [\sum_{i=1}^{n} B_i>C]\sim \Pr [\max (B_1,B_2,\dots
,B_n)>C]$ as $C\rightarrow \infty$, implying asymptotic relation
\begin{equation*}
\Pr \left[\sum_{i=1}^{n}B_i>C,B_{i}\le C\text{ for every $1\le i\le
n$}\right]=o(\Pr [B>C])\;\;\text{as}\;\;\text{$C\rightarrow
\infty$.}
\end{equation*}
In order to show that $n$ can be replaced by $N_{s,n}$ in the
above inequality, we need to integrate it with respect to the
density of $N_{s,n}$, i.e.,
\begin{align*}
&\Pr \left[\sum_{i\in {\cal N}_{s,n}\cup \{n\}}B_i>C,B_{i}\le
C\text{ for every
$i\in {\cal N}_{s,n}\cup \{n\}$}\right]\\
&\qquad\qquad\qquad=\sum_{k=0}^{\infty}\Pr [N_{s,n}=k]\Pr
\left[\sum_{i=1}^{k+1}B_i>C,B_{i}\le C\text{ for every $i=1,\dots
,k+1$}\right].
\end{align*}
Note that on the left hand side of the previous equation index $i$
can take negative values. Next, due to the lemma stated by Kesten
(see Lemma $7$,  pp.$149$ of \cite{ATN72}), for any $\epsilon>0$
there exists a positive constant $K(\epsilon)$ such that
\begin{eqnarray*}
\frac{\Pr[\sum_{i=1}^{k}B_i>C,B_{i}\le C\text{ for every $1\le
i\le k$}]}{\Pr [B>C]} \le \frac{\Pr[\sum_{i=1}^{k}B_i>C]}{\Pr
[B>C]} \le K(\epsilon)(1+\epsilon)^k,
\end{eqnarray*}
for any integer $k$ and all capacity values $C<\infty$. Then,
since the probability generating function $\expect z^{N_{s,n}}$ is
finite for any $z\in \mathbb{C}$ (see Theorem 1 in \cite{TAK80}
and Theorem 5 in \cite{LKT90} for the detailed proof), we have
$\sum_{k=0}^{\infty} \Pr [N_{s,n} =k](1+\epsilon)^{k} <\infty$.
Therefore, by applying the dominated convergence theorem, we
conclude that
\begin{align}
&\lim_{C\rightarrow \infty}\frac{\Pr \left [\sum_{i\in
N_{s,n}}B_i+B_n>C,B_{i}\le C\text{ for every $i\in
{\cal N}_{s,n}\cup \{n\}$}\right]}{\Pr [B>C]}\nonumber\\
&=\lim_{C\rightarrow \infty}\sum_{k=0}^{\infty}\frac{\Pr
[N_{s,n}=k]\Pr \left [\sum_{i=1}^{k+1}B_i
>C,B_i\le C\text{ for every $1\le i\le
k+1$}\right]}{\Pr [B>C]}\nonumber\\
&=0,\label{eq:ref2}
\end{align}
which in conjunction with (\ref{eq:split1}) and (\ref{eq:lowbd5}),
completes the proof of this theorem. \qed

\noindent{\bf Remark:} It may appear surprising that the
performance of the loss network from above does not depend on
engagement durations, as long as they have finite mean. In
addition, the result is quite general and provides the asymptotic
result for a large (subexponential) class of possible resource
requirement distributions.

\subsection{Advance reservations}
\label{sec:advancereservation}

Using the result of Theorem \ref{theorem:th1} and observations
from the previous remark, we extend the loss networks model to
allow requests to become effective with some delay with respect to
the moments of their arrivals. In particular, a request that
arrives at time $\tau_n$ and requires $B_n$ amount of resources
for some random time $\theta_n$ starting from the moment
$\tau_n+D_n$ is accepted if previously admitted resource
requirements allow that; otherwise, it is rejected. In other
words, a request arriving at $\tau_n$ is lost if at any moment of
time in interval $(\tau_n+D_n, \tau_n+D_n+\theta_n)$ the total
amount of active requirements requested prior to $\tau_n$ exceeds
$C-B_n$. First, note that $B_n
>C$ implies the loss of $n$th request and, therefore, it is
straightforward to conclude that the blocking probability in the
system with advance reservations can be lower bounded by $\Pr
[B>C]$.

Next, we discuss the idea behind proving the upper bound on the
blocking probabilities. By applying sample path arguments one can
show that, at any moment of time, the amount of active resources
in the previously described system with advance reservations can
be bounded from above by the amount of active resources in another
system of unlimited capacity, without advance reservations, with
resource holding times $D_n+\theta_n$ for every $n$, and with
requests for resources being sampled from the original process
$\{B_n\}$ whenever the corresponding requirements are less or
equal to $C$. Equivalently, the blocking probability in the system
with advance reservations can be bounded from above by
\begin{equation*}
\Pr \left [\sum_{i\in {\cal N}^{(C)}_{s,n} (\theta +D)} B_i +B_n
>C \right ],
\end{equation*}
where ${\cal N}^{(C)}_{s,n} (\theta+D)$ is a set of request
indices $i<n$ that are active at time $\tau_n$, whose requirements
are less or equal to $C$ and holding times last throughout the
interval $(\tau_i,\tau_i+D_i+\theta_i)$, assuming that there is an
unlimited resource capacity.

Finally, by using the previous discussion, the properties of
$\{B_n\}$, $\{\theta_n\}$ and $\{\tau_n\}$ as introduced at the
beginning of this section, assuming that reservation times
$\{D_n\}$, $\expect D_n <\infty$, are i.i.d. and independent from
$\{B_n\}$, $\{\theta_n\}$ and $\{\tau_n\}$, and applying the
identical arguments as in the proof of Theorem \ref{theorem:th1},
we obtain the following result:
\begin{corollary}
The blocking probability in the system with advance reservations
approaches $\Pr [B>C]$ as $C\rightarrow \infty$.
\end{corollary}

\section{Acquiring resources of different types (loss networks case)}
\label{sec:network}

Assume that there are $K\in \nat$ resource types with capacities
$C_1,\dots ,C_K$. Again, requests arrive at
$\{\tau_n,-\infty<n<\infty\}$, which represent a renewal process
with rate $0<\lambda =1/\expect [\tau_1 -\tau_0] <\infty$. There
are $M<\infty$ request types and, given an arrival, the request is
of type $l$, $1\le l\le M$, with probability $p_l$,
$p_1+\dots+p_M=1$, independent from $\{\tau_n\}$. We will use
random variables $J_n\in \{1,2,\dots ,M\}$ to denote the type of
the request arriving at $\tau_n$. Furthermore, let
$B_n^{(J_n,1)},\dots ,B_n^{(J_n,K)}$ represent amounts of required
resources of each type at time $\tau_n$ and let
$\theta_n^{(J_n)}$, $\expect \theta_n^{(J_n)}<\infty$, be the
corresponding random duration. We assume that sequences
$\{(B_n^{(J_n,1)},\dots , B_n^{(J_n,K)})\} ,\{\theta_n^{(J_n)}\}$
are mutually independent and independent from $\{\tau_n\}$. Given
the event $\{J_n=l\}$, resource requirements $B_n^{(l,i)}$, $1\le
i\le K$, are mutually independent nonnegative random variables
drawn from distributions $F_{l,i}$, $1\le i\le K$; if a request
does not require resources of type $i$ then $B_n^{(l,i)}=0$ a.s.,
$-\infty <n<\infty$. Only if there is enough capacity available,
the request arriving at time $\tau_n$ will be accepted and all of
the engaged resources will be occupied for the duration of
$\th_n^{(J_n)}$; otherwise, the request is rejected.

Our goal is to estimate the blocking probability in a system
described above. Define $Q^{(1)}_n,\dots ,Q_n^{(K)}$ to be amounts
of resources of each type that a request arriving at time $\tau_n$
finds engaged. Note that $Q^{(i)}_n$, $1\le i\le K$, are mutually
dependent and, as pointed out in the Introduction, it is hard to
compute the blocking probability of this system explicitly. Using
analogous arguments as in the case of a single resource type (see
the Appendix), one can show that the stationary distribution of
$Q^{(i)}_n$, $1\le i\le K$, exists. Probability that the request
arriving at time $\tau_n$ is blocked equals to
\begin{equation} \label{eq:lossmult} \Pr [\cup_{1\le i\le K}
\{Q_n^{(i)}+B_n^{(J_n,i)}>C_i\}],
\end{equation}
and our goal again is to estimate its value as $\min_i C_i$ grows
large.

Asymptotic estimates derived in this section hold under the
following assumption: \vspace{0.2in}

\noindent{\bf Assumptions:} For each resource type $1\le i\le K$,
let ${\cal L}_i$ and ${\cal H}_i$ be two disjoint sets of request
types ($|{\cal L}_i\cup {\cal H}_i|=M$) satisfying:
\begin{itemize}
\item  Assume that there exists at least one resource type that is
accessed by subexponentially distributed resource requirements,
which implies $|{\cal H}_i|>0$ for some $1\le i\le K$; \item For
every $l\in {\cal H}_i \neq \emptyset$, there exists a
subexponential distribution $F_i\in {\cal S}$ such that $\Bar
{F}_{l,i}(x)\sim c_{l,i}\Bar {F}_i(x)$ as $x\rightarrow \infty$
with $c_{l,i}>0$; \item There exists a subexponential random
variable $L\in {\cal S}$ that satisfies
\begin{eqnarray*}
\Pr[L>x]\ge \max_{1\le i\le K,l\in {\cal L}_i}  \Pr
[B_n^{(J_n,i)}>x|J_n=l]\;\;\text{for all $x>0$,}
\end{eqnarray*}
and $\Pr [L>x]=o(\Bar {F}_i (x))$ as $x\rightarrow \infty$ for all
$i\in \{j|{\cal H}_j\neq \emptyset\}$.
\end{itemize}

\noindent {\bf Remark: } In the preceding assumptions, we require
the resource requirement distributions to be asymptotically
comparable. For each $1\le i\le K$, ${\cal H}_i$ contains tail
dominant subexponential distributions that are asymptotically
proportional to each other. On the other hand, the only assumption
imposed on the distributions in ${\cal L}_i$, $1\le i\le K$, is
that there is a subexponential tail that asymptotically dominates
them. This asymptotic tail comparability is necessary for our main
result to hold. In particular, these conditions are extensively
used in (\ref{eq:chtermth2}) - (\ref{eq:unifbound}) of the proof
of Theorem \ref{theorem:th2}.

Next, we prove the following lemma that investigates summations of
random variables with different tail distributions.
\begin{lemma}
\label{lem:tail} Suppose that $X_i, 1\le i\le n$, are independent
random variables with corresponding tail distributions ${\bar
F}_i(x)$, $1\le i\le n$. If there exists $F\in {\cal S}$ such that
${\bar F}_i(x)\sim c_i {\bar F}(x)$ as $x\rightarrow \infty$ with
$c_i\ge 0, 1\le i\le n$, and $\sum_{i=1}^n c_i>0$, then the
following asymptotic relation holds:
\begin{equation}
\label{eq:mainlemma} \Pr \left [ \sum_{i=1}^n X_i
>x, X_i \le x, 1\le i\le n \right ]= o({\bar F}(x)
)\;\;\text{as}\;\;\text{$x\rightarrow \infty$}.
\end{equation}
\end{lemma}
\proof Note that
\begin{equation*}
\Pr \left [\sum_{i=1}^{n} X_i>x\right ]=\Pr \left [\sum_{i=1}^{n}
X_i>x, X_i\le x, 1\le i\le n\right ]+\Pr \left [\sum_{i=1}^{n}X_i
>x,\cup_{i=1}^{n} \{X_i>x\}\right ].
\end{equation*}
Then, the previous expression, $\cup_{i=1}^{n}
\{X_i>x\}\subset\{\sum_{i=1}^{n}X_i
>x\}$, independence of $X_i$s, as well as Lemmas 4.2 and 4.5 of
\cite{AHK94},
imply (\ref{eq:mainlemma}). \qed

First, we estimate the asymptotic lower bound for the expression
in (\ref{eq:lossmult}). By using our model assumptions,
$\{B_n^{(J_n,i)}>C_i\}\subset \{Q_n^{(i)}+B_n^{(J_n,i)}>C_i\}$ and
independence, we obtain
\begin{equation}
\label{eq:finallow} \Pr [\cup_{1\le i\le K}
\{Q_n^{(i)}+B_n^{(J_n,i)}>C_i\}]\ge \Pr [\cup_{1\le i\le K}
\{B_n^{(J_n,i)}>C_i\}]\sim \sum_{i=1}^{K} \sum_{l\in {\cal H}_i}
p_l \Bar {F}_{l,i} (C_i),
\end{equation}
as $\min_i C_i\rightarrow \infty$.

Next, we estimate the asymptotic upper bound for the expression in
(\ref{eq:lossmult}). Using the union bound yields
\begin{align}
\Pr [\cup_{1\le i\le K} \{Q_n^{(i)}+B_n^{(J_n,i)}>C_i\}] &\le
\sum_{i=1}^{K} \Pr
[Q_n^{(i)}+B_n^{(J_n,i)}>C_i].\label{eq:crucialbd}
\end{align}
Similarly as in (\ref{eq:refth11}) of Theorem~\ref{theorem:th1},
for each resource $1\le i\le K$,
\begin{equation}
\label{eq:sampleth2} \Pr [Q_n^{(i)}+B_n^{(J_n,i)}>C_i]\le \Pr
\left [\sum_{l\in {\cal L}_i}\sum_{j\in {\cal N}_{s,n}^{(l,C_i)}}
B_j^{(l,i)}+ \sum_{l\in {\cal H}_i} \sum_{j\in {\cal
N}^{(l,C_i)}_{s,n}}B_j^{(l,i)}+B_n^{(J_n,i)}>C_i\right ],
\end{equation}
where ${\cal N}_{s,n}^{(l,C_i)}$, $1\le l\le M$, are sets of
indices $j<n$ defined as
\begin{equation*}
{\cal N}_{s,n}^{(l,C_i)}\eqdef \{ j<n |J_j=l,B_j^{(l,i)}\le C_i,
\theta_j^{(l)}>\tau_n -\tau_j\}.
\end{equation*}
In the previous expressions we bounded the amount of allocated
resources that are active at time $\tau_n$ by the corresponding
quantity in another system of infinite capacity where the
corresponding request process is sampled from the original
$\{B^{(J_n,i)}_n\}$, $1\le i\le K$, whenever the corresponding
requirements are less than or equal to $C_i$, $1\le i\le K$.

In the rest of the proof, we derive an asymptotic estimate for the
expression in (\ref{eq:sampleth2}). After conditioning on $\{
N_{s,n}^{(1,C_i)}=n_1,\dots ,N_{s,n}^{(M,C_i)}=n_M\}$
($N_{s,n}^{(l,C_i)}\eqdef |{\cal N}_{s,n}^{(l,C_i)}|$, $1\le l\le
M$), we obtain
\begin{align}
&\Pr [Q_n^{(i)}+B_n^{(J_n,i)}>C_i]\nonumber\\
&\le \sum_{0\le n_1,\dots ,n_M<\infty}\Pr
[N_{s,n}^{(1,C_i)}=n_1,\dots ,N_{s,n}^{(M,C_i)}=n_M]\nonumber\\
&\;\;\times \Pr \left [\sum_{l\in {\cal H}_i}\sum_{j=1}^{n_l}
B_{(j)}^{(l,i)} + \sum_{l\in {\cal L}_i}\sum_{j=1}^{n_l}
B_{(j)}^{(l,i)}+ B_n^{(J_n,i)} >C_i,\text{ $B_{(j)}^{(l,i)}\le
C_i$, $1\le j\le n_l$, $1\le l\le M$}\right
],\label{eq:mainth2term}
\end{align}
where $B^{(l,i)}_{(j)}\eqd B_k^{(l,i)}$, $k\in {\cal
N}^{(l,C_i)}_{s,n}$, $j=1,\dots ,n_l$, are independent replicas of
requests in ${\cal N}^{(l,C_i)}_{s,n}$. Next, after conditioning
on $\{J_n=m\}$, $m=1,\dots ,M$, and then on $B_n^{(m,i)}$ being
smaller or larger than $C_i$, we can further upper bound the
conditional blocking probability in (\ref{eq:mainth2term}) as
\begin{align}
&\Pr \left [\sum_{l\in {\cal H}_i}\sum_{j=1}^{n_l} B_{(j)}^{(l,i)}
+ \sum_{l\in {\cal L}_i}\sum_{j=1}^{n_l} B_{(j)}^{(l,i)}+
B_n^{(J_n,i)} >C_i,\text{ $B_{(j)}^{(l,i)}\le C_i$, $1\le j\le
n_l$, $1\le l\le M$}\right ]\le \nonumber\\
&\sum_{m=1}^{M} p_m\Pr \left [\sum_{l\in {\cal
H}_i}\sum_{j=1}^{n_l} B_{(j)}^{(l,i)} + \sum_{l\in {\cal
L}_i}\sum_{j=1}^{n_l} B_{(j)}^{(l,i)}+ B_n^{(m,i)} >C_i,\text{
$B_{(j)}^{(l,i)}\le C_i$, $1\le j\le n_l$, $1\le l\le
M$, $B_n^{(m,i)}\le C_i$}\right ]\nonumber\\
&\qquad\qquad + \sum_{m=1}^{M} p_m\Pr [B_n^{(m,i)}
>C_i].\label{eq:mainsplitth2}
\end{align}
Thus, the probabilities in the first term on the right hand side
of the previous expression can be expressed in the form
\begin{equation}
\label{eq:chtermth2} \Pr \left [\sum_{l\in {\cal
H}_i}\sum_{j=1}^{n'_l} B_{(j)}^{(l,i)} + \sum_{l\in {\cal
L}_i}\sum_{j=1}^{n'_l} B_{(j)}^{(l,i)}>C_i,\text{
$B_{(j)}^{(l,i)}\le C_i$, $1\le j\le n'_l$, $1\le l\le M$}\right
],
\end{equation}
where $n'_l=n_l$ for $l\not =m$ and $n'_l=n_l+1$ for $l=m$.

Next, in order to estimate the asymptotic upper bound of the term
in (\ref{eq:chtermth2}), Assumptions enable us to distinguish
between two cases: (i) ${\cal H}_i =\emptyset$ or $\sum_{l\in
{\cal H}_i}n'_l =0$, and (ii) ${\cal H}_i \neq \emptyset$ and
$\sum_{l\in {\cal H}_i}n'_l
>0$.

\noindent (i): If ${\cal H}_i =\emptyset$ or $\sum_{l\in {\cal
H}_i}n'_l =0$, we have that probability in (\ref{eq:mainsplitth2})
can be upper bounded as
\begin{equation*}
\Pr \left [ \sum_{l\in {\cal L}_i}\sum_{j=1}^{n'_l}
B_{(j)}^{(l,i)}>C_i\right ]\le \Pr \left [ \sum_{l\in {\cal
L}_i}\sum_{j=1}^{n'_l} L^{(l,i)}_{(j)}>C_i\right ],
\end{equation*} where in the inequality above we used Assumptions and
introduced $L_{(j)}^{(l,i)}$ to be independent r.v.s equal in
distribution to $L$. Hence, since $L^{(l,i)}_{(j)}$ are
subexponential, we obtain
\begin{equation}
\lim_{C_i\rightarrow \infty}\frac{\Pr \left [ \sum_{l\in {\cal
L}_i}\sum_{j=1}^{n'_l} B_{(j)}^{(l,i)}>C_i\right ]}{\Pr
[L>C_i]}\le \sum_{l\in {\cal L}_i} n'_l.\label{eq:lightpieceth2}
\end{equation}

\noindent (ii): If ${\cal H}_i \neq \emptyset$ and $\sum_{l\in
{\cal H}_i}n'_l >0$, using Assumptions and Lemma 1, we derive the
following asymptotic upper bound
\begin{equation}
\Pr \left [\sum_{l\in {\cal H}_i}\sum_{j=1}^{n'_l} B_{(j)}^{(l,i)}
+ \sum_{l\in {\cal L}_i}\sum_{j=1}^{n'_l}
B_{(j)}^{(l,i)}>C_i,\text{ $B_{(j)}^{(l,i)}\le C_i$, $1\le j\le
n'_l$, $1\le l\le M$}\right ]=o({\bar F}_i(C_i)),
\label{eq:estth21}
\end{equation}
as $C_i\rightarrow \infty$.

Thus, in (\ref{eq:chtermth2})-(\ref{eq:estth21}) we obtained upper
bounds and their asymptotic estimates for the conditional blocking
probabilities in the first term of (\ref{eq:mainsplitth2}) that
hold for any finite nonnegative integers $n_1,\dots,n_M$. Thus, in
view of (\ref{eq:mainth2term}), in order to estimate an asymptotic
upper bound of $\Pr [Q_n^{(i)}+B_n^{(J_n,i)}>C_i]$, we need to
integrate probabilities in (\ref{eq:chtermth2}) with respect to
densities of r.v.s $N^{(l,C_i)}_{s,n}$, $l=1,\dots ,M$. In this
regard, note that in the case where ${\cal H}_i\not =\emptyset$,
by Assumptions, the term in (\ref{eq:chtermth2}) can be upper
bounded as
\begin{align}
&\Pr \left [\sum_{l\in {\cal H}_i}\sum_{j=1}^{n'_l}
B_{(j)}^{(l,i)} + \sum_{l\in {\cal L}_i}\sum_{j=1}^{n'_l}
B_{(j)}^{(l,i)}>C_i,\text{ $B_{(j)}^{(l,i)}\le C_i$, $1\le j\le
n'_l$, $1\le l\le M$}\right ]\nonumber\\
&\le \Pr \left [\sum_{l\in {\cal H}_i}\sum_{j=1}^{n'_l}
B_{(j)}^{(l,i)} + \sum_{l\in {\cal L}_i}\sum_{j=1}^{n'_l}
L^{(l,i)}_{(j)}>C_i\right ],\label{eq:upbd100th2}
\end{align}
where, as before, $L^{(l,i)}_{(j)}$ are independent r.v.s equal in
distribution to $L$. Furthermore, since $\Pr [L>x]=o({\Bar F}_i
(x))$ as $x\rightarrow \infty$, there exists a large enough finite
integer $H$ such that $\Pr [L>x]\le H{\Bar F}_i (x)$ for all $x\ge
0$. Therefore, for any $x\ge 0$, one can write
\begin{equation}
\label{eq:stbdth2} \Pr [L>x]\le H{\Bar F}_i (x)=\Pr \left
[\cup_{1\le r\le H}\{{\Hat B}^{(i)}_r >x\}\right ]\le \Pr \left
[\sum_{r=1}^{H} {\Hat B}^{(i)}_r >x\right ],
\end{equation}
where ${\Hat B}^{(i)}_r$, $1\le r\le H$, are independent r.v.s
having cumulative distribution function $F_i$. Now, in view of
(\ref{eq:stbdth2}), each of random variables $L^{(l,i)}_{(j)}$ in
(\ref{eq:upbd100th2}) can be stochastically upper bounded by a
random variable that is equal in distribution to $\sum_{r=1}^{H}
{\Hat B}^{(i)}_r$. Thus, if we introduce $Y_j$, $j\ge 1$, to be
independent r.v.s equal in distribution to $\sum_{r=1}^{H} {\Hat
B}^{(i)}_r$, we obtain
\begin{equation*}
\Pr \left [\sum_{l\in {\cal H}_i}\sum_{j=1}^{n'_l} B_{(j)}^{(l,i)}
+ \sum_{l\in {\cal L}_i}\sum_{j=1}^{n'_l}
L_{(j)}^{(l,i)}>C_i\right ]\le \Pr \left [\sum_{l\in {\cal
H}_i}\sum_{j=1}^{n'_l} B_{(j)}^{(l,i)}+\sum_{j=1}^{\sum_{l\in
{\cal L}_i} n'_l} Y_j>C_i\right ],
\end{equation*}
which in conjunction with point (b) of Lemma 4.2 in \cite{AHK94}
implies that for any $\epsilon >0$ there exist a finite constant
$K_{\epsilon}$ such that
\begin{align}
\Pr \left [\sum_{l\in {\cal H}_i}\sum_{j=1}^{n'_l} B_{(j)}^{(l,i)}
+ \sum_{l\in {\cal L}_i}\sum_{j=1}^{n'_l}
L_{(j)}^{(l,i)}>C_i\right ]&\le \Pr \left [\sum_{l\in {\cal
H}_i}\sum_{j=1}^{n'_l} B_{(j)}^{(l,i)}+\sum_{j=1}^{\sum_{l\in
{\cal L}_i}n'_l}
Y_j>C_i\right ]\nonumber\\
&\le K_{\epsilon} (1+\epsilon)^{\sum_{l\in {\cal H}_i} n'_l
+\sum_{l\in {\cal L}_i} n'_l} {\Bar F}_i(C_i),\label{eq:unifbound}
\end{align}
for any $C_i<\infty$. Similarly, in cases where ${\cal H}_i
=\emptyset$, we could apply the stochastic dominance
$B_{(j)}^{(l,i)}\leqd L^{(l,i)}_{(j)}$, $l\in {\cal L}_i$, where
$L^{(l,i)}_{(j)}$ are, as before, independent subexponential
random variables equal in distribution to $L$. Then, by Kesten's
lemma (see Lemma 7 on page 149 of \cite{ATN72}), the analogous
bound to the one in (\ref{eq:unifbound}) follows.

Finally, since (\ref{eq:unifbound}) bounds uniformly probabilities
in (\ref{eq:chtermth2}) for all $C_i<\infty$ and $n'_l$, $1\le
l\le M$, in conjunction with (\ref{eq:mainsplitth2}),
(\ref{eq:mainth2term}), $N_{s,n}^{(l,C_i)}\le N_{n}^{(l,\infty)}$
a.s. and existence of $ \expect z^{N_{n}^{(l,\infty)}}$ for all
$z\in \mathbb{C}$, $1\le l\le M$, (see Theorem 1 in \cite{TAK80}
and Theorem 5 in \cite{LKT90}), one can apply the dominated
convergence theorem and conclude
\begin{equation*} \lim_{C_i\rightarrow \infty}\frac{\Pr [Q_n^{(i)}+B_n^{(J_n,i)}>C_i]}{\sum_{l\in {\cal H}_i} p_l \Bar {F}_{l,i} (C_i)}\le 1[ {\cal H}_i \neq \emptyset].
\end{equation*}
Next, by adding asymptotic estimates for all $1\le i\le K$, in
conjunction with (\ref{eq:finallow}), we complete the proof of the
following result:
\begin{theorem}
\label{theorem:th2} For the request model introduced in this
section, under the conditions imposed by Assumptions, the
stationary blocking probability for general loss networks
satisfies
\begin{equation*}
\Pr [\cup_{1\le i\le K} \{Q_n^{(i)}+B_n^{(J_n,i)}>C_i\}]\sim
\sum_{i=1}^{K}\sum_{l\in {\cal H}_i}p_lc_{l,i} \Bar
{F}_{i}(C_i)\;\;\text{as}\;\;\text{$\min_i C_i\rightarrow \infty$.}
\end{equation*}
\end{theorem}

\section{Numerical examples}
\label{sec:experiment}

In this section, with two simulation experiments, we demonstrate the
accuracy of our asymptotic formulas, proved in Theorems
\ref{theorem:th1} and \ref{theorem:th2}. Our goal is to show that
even though our results are asymptotic, the derived estimates match
experiments with high accuracy even for systems with finite support
demand distributions and moderately large capacities.

\begin{figure}[htb]
\epsfxsize=8cm \centerline { \epsffile{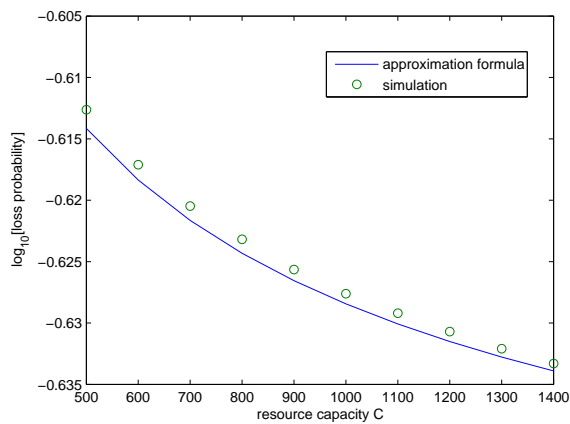}}
\caption{Illustration for Example 1} \label{fig:example1}
\end{figure}

\begin{figure}[htb]
\epsfxsize=8cm \centerline { \epsffile{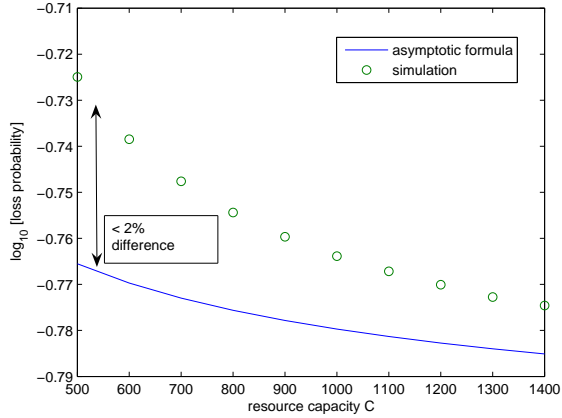}}
\caption{Illustration for Example 2} \label{fig:example2}
\end{figure}

In each experiment, in order for the system to reach stationarity,
we let the first $10^{8}$ arrivals to be a warm-up time. By
repeating many experiments, we observe that longer warm-up times
do not lead to improved results.
Then, we count the number of blocked requests among next $10^{9}$
arrivals. In both of the experiments below, measurements are
conducted for capacities $C=500+100j,0\le j\le 9$, where the
starting value of $C=500$ is set to be slightly larger than the
effective systems load $\lambda \expect[\theta_n] \expect [B_n]$.
Simulation results are presented by symbol ``o'' in
Figures~\ref{fig:example1} and \ref{fig:example2}, while our
approximations, estimates obtained in Theorems \ref{theorem:th1}
and \ref{theorem:th2}, are the solid lines on the same figures.
Note that in order to emphasize the difference and to observe a
range of blocking probabilities we are trying to estimate, we
present base $10$ logarithm of the obtained values.

\noindent{\bf Example 1} Consider the case of a single resource
type of capacity $C$. Let requests for resources arrive at Poisson
time points with rate $\lambda = 1$. In addition, we assume that
engagement durations are exponentially distributed with mean
$1/\mu =1$. Next, let request requirements $B_n$ be drawn from a
finite support distribution, where $\Pr [B_n
=i]=\frac{0.3}{i^{1.5}}$, $1\le i\le 1999$, and $\Pr
[B_n=2000]=1-\Pr [B_n <2000]$ (power law distribution). Effective
load in this example is $\lambda \expect [\theta_n] \expect
[B_n]\approx 485.8$. Experimental results are presented in Figure
\ref{fig:example1}. Even though we start measuring rejections at
capacities that are slightly larger than the mean requirement
value, our approximation $\Pr[B_n>C]$  is very close to
experimental results. In particular, the relative approximation
error is less than $1\%$ for $C=500$, and for capacity values
larger or equal to $C=1400$ this error is less than $0.3\%$.

\noindent{\bf Example 2} In this example, we consider the case of
two resource and two request types. Furthermore, we assume that
resource capacities are the same $C=C_1=C_2$. The frequencies of
requests of types 1 and 2 are $p_1=0.3$ and $p_2=0.7$
respectively. Assume that the arrival points are separated by a
fixed, unit length of time, i.e., $\tau_n-\tau_{n-1}=1$ for all
$n$. Type $1$ request durations satisfy $\theta^{(1)}_i\sim
exp(4)$ and type 2 request holding times are drawn from the
uniform distribution on $[0,40]$, i.e., $\theta^{(2)}_i\sim
Unif([0,40])$. Resource requirements corresponding to engagements
of type 1 are distributed as $\Pr [B^{1,1}=1]=0.8$, $\Pr
[B^{1,1}=i]=0.15e^{-\sqrt {i}}$, $2\le i\le 1999$ and $\Pr
[B^{1,1}=2000]=1-\Pr [B^{1,1}<2000]$ for the type 1 resources, and
$\Pr [B^{1,2}=50]=1$ for type 2 resources. Requests of type 2
require resources according to $\Pr
[B^{2,1}=i]=geom^{i-1}(1-geom)$, $1\le i\le 1999$, $\Pr
[B^{2,1}=2000]=1-\Pr [B^{2,1}<2000]$, where $geom=0.6$ for
resources of type 1, and $\Pr [B^{2,2} =i]=\frac{0.3}{i^{1.5}}$,
$1\le i\le 1999$, $\Pr [B^{2,2}=2000]=1-\Pr [B^{2,2} <2000]$ for
type 2 resources. Our asymptotic results suggest that the blocking
probability should be characterized by the heaviest tailed demand
distributions. The results of this experiment are presented in
Figure \ref{fig:example2}. As in the previous case, we obtain a
very accurate agreement between our approximation and the
simulation. The relative approximation error in this case does not
exceed $2\%$ and is getting smaller as resource capacities grow.

\noindent{\bf Remark:} (i) We would like to point out that the
accuracy of experimental results directly depends on the
approximation errors (\ref{eq:refth11}) and (\ref{eq:estth21}),
depending on the simulated scenarios. These errors highly depend
on the tail properties of the resource requirements distributions.
More specifically, under fairly general assumptions, the heavier
the dominant tail of the resource requirement distribution is, the
smaller would be the relative approximation error. For detailed
explanations, a reader is referred to Section 1.3.2 of
\cite{EKM97}. (ii) Note that our main results estimate the
stationary blocking probability and, as we commented earlier, are
indifferent to distributional properties of holding times. For
that reason, as long as one can claim that the measurements are
conducted in stationarity, the transience should not affect
experimental results.

\section{Concluding remarks}
\label{sec:concl}

In this paper, we consider loss networks with reusable resources
and finite resource capacities and estimate the probability that a
request is rejected due to insufficient amount of resources at
points of their arrivals. Assuming a renewal process of request
arrivals, subexponential resource requirements and generally
distributed activity durations, we show that the asymptotic
blocking probability for a wide class of analyzed systems can be
fully estimated using resource requirement distribution,
independent from other system's properties. In particular, we show
that the blocking probability behaves as the asymptotically
dominant tail of the resource requirement distribution.

The model we study can be applied to a wide range of applications.
Historically, loss networks (in particular, Erlang loss networks)
are widely used for modeling communication networks. Later,
through the development of new services applications such as
workforce management with similar modeling properties, the
importance of accurately estimating blocking probabilities of
general loss networks has become significant. In this regard, we
investigate loss networks with various request types and possibly
highly variable random amounts of required resources. In addition,
we research the possibility of incorporating random advance
reservations for incoming requests. These results should be of
great interest to an emerging research community. Although our
results are intended mainly for qualitative purposes, numerical
examples demonstrate an excellent match between derived formulas
and simulated systems performance, hence strongly suggesting their
application.

\section*{Appendix}

In this section we prove the existence of the stationary blocking
probabilities in (\ref{eq:deflossrate}). Using the model
description from Section \ref{sec:model}, we observe the system at
the moments of request arrivals. Then, we define a discrete time
process $X_n\triangleq (N^{(C)}_{n}, B_i, E_i, i=1, 2, \cdots,
N^{(C)}_{0,n})$, where $E_i$s represent times that elapsed in
processing requests in the system by time $\tau_n$; furthermore,
$N^{(C)}_{0,n}$ is the number of active requests at the moment of
$n$th arrival. Note that $\{X_n\}_{n\ge 0}$ is a discrete time
Markov chain with state space $\Omega\eqdef \nat_{0}\times
\ell_\infty \times \ell_\infty$, where $\ell_\infty$ denotes the
Banach space of the infinite sequence of real numbers equipped
with the supreme norm; let $\omega_0\in\Omega$ denote the state
with no active requests. We start observing the system at the
moment $\tau_0$ of $0$th arrival and denote the initial state by
$X_0=(N^{0},B_i^{0},E_i^{0},i=1,\dots,N^{0})$ drawn from some
arbitrary distribution $P_0$, where $\expect B_i^{0}<\infty$,
$\expect \theta^{0}_i<\infty$. Next, define ${\mathcal F}$ to be
the Borel field of $\Omega$, and let $P_n (x_0,A)$, $x_0\in
\Omega$, $A\in {\mathcal F}$, represent a transition probability
of the Markov chain $X_n$ into set $A$ in time $n$, starting from
state $x_0$. Let $P_n$ be the probability distribution of $X_n$.

Now, in order to prove the existence of a unique stationary
distribution for the Markov chain $\{X_n\}$, we use a discrete
version of Theorem 1 in \cite{SEV57}, which we state next for
reasons of completeness.

\noindent{\bf Theorem:} A Markov chain homogeneous in time has a
unique stationary distribution which is ergodic if, for any
$\epsilon >0$, there exists a measurable set $S$, a probability
distribution $R$ in $\Omega$, and $n_1>$, $k>0$, $K>0$ such that
\begin{itemize}
\item $k R(A)\le P_{n_1} (x,A)$ for all points $x\in S$ and
measurable sets $A\subset S$; for any initial distribution $P_0$
there exists $n_0$ such that for any $n\ge n_0$ \item $P_n (S)\ge
1-\epsilon$, \item $P_n (A) \le K R(A) +\epsilon$ for all measurable
sets $A\subset S$.
\end{itemize}

\proof The proof follows identical arguments as in \cite{SEV57}
translated into discrete setting.

Next, we need to show that Theorem 1 holds for the process
investigated in this paper; in particular, we will consider a
common resource pool case. The proof follows the similar reasoning
as in Theorems 4 and 5 of \cite{SEV57}.

Define set $S(\psi, \beta, \delta)$ as
\begin{equation*}
S(\psi,\beta, \delta)\eqdef  \{N^{(C)}_{0,n}\le \psi, 0\le B_i\le
\beta, 0\le E_i\le \delta, i\in {\cal N}^{(C)}_{0,n}\}
\end{equation*}
for some positive finite constants $\psi, \beta, \delta$.

Now, we show that for any $\epsilon >0$, there exists $S(\psi,
\beta,\delta) \subset \Omega$ such that for any initial distribution
$P_0$ there exists $n_0$ such that for all $n\ge n_0$
\begin{equation}
\label{eq:point2} P_n (S(\psi,\beta,\delta))\ge 1-\epsilon.
\end{equation}
Note that
\begin{align}
&P_n (\Bar {S} (\psi,\beta,\delta))\le \Pr [\cup_{i\in {\cal
N}^{(C)}_{0,n}} \{\theta_i > \delta \}, N^{(C)}_{0,n}\le \psi] + \Pr
[\cup_{i\in {\cal N}^{(C)}_{0,n}} \{B_i > \beta\},N^{(C)}_{0,n}\le
\psi]+\Pr [N^{(C)}_{0,n}>\psi]\nonumber\\
&\le \psi \Pr [\theta_i > \delta]+\psi \Pr [B_i >\beta] + \Pr
[N^{(C)}_{a,n}+N^{0}_{0,n}>\psi],\label{eq:point2p1}
\end{align}
where $N^{(C)}_{a,n}$ represents the number of active requests at
$\tau_n$ that originated from $n$ arrivals at $\tau_0,\dots
\tau_{n-1}$, and the rest of active requests at $\tau_n$,
$N^{0}_{0,n}=N^{(C)}_{0,n}-N^{(C)}_{a,n}$ are those that were active
at the initial point $\tau_0$ and are still processed at the moment
of $n$th arrival. Next, since
\begin{align}
&\Pr [N^{(C)}_{a,n}+N^{0}_{0,n}>\psi]\le \Pr \left
[N^{(C)}_{a,n}>\frac{\psi}{2}\right ] +\Pr \left [N^{0}_{0,n}
>\frac{\psi}{2}\right ]\nonumber\\
&\le \Pr \left [N^{(\infty)}_n>\frac{\psi}{2}\right ]+\Pr
\left [\sum_{i=1}^{N^{0}}1[\theta_i^{0}>\tau_n-\tau_0]>\frac{\psi}{2}\right ]\nonumber\\
&\le \Pr \left [N^{(\infty)}_n>\frac{\psi}{2}\right
]+\sum_{m=0}^{\infty}\Pr [N^{0}=m]\Pr \left
[\sum_{i=1}^{m}1[\theta_i^{0}>(1-\epsilon_1)n\expect
[\tau_1-\tau_0]]>\frac{\psi}{2}\right ]\nonumber\\
&\qquad\qquad+\Pr [\tau_n-\tau_0< (1-\epsilon_1)n\expect
[\tau_1-\tau_0]],\label{eq:initstat1}
\end{align}
where in the previous inequalities $0<\epsilon_1 <1$ is an arbitrary
constant and we used $N^{(\infty)}_n \ge N^{(C)}_{a,n}$ a.s., where
$N^{(\infty)}_n$ is defined as in the proof of Theorem
\ref{theorem:th1}.

Now, we prove that there exists $\psi=\psi_0$ large enough such that
(\ref{eq:initstat1}) is bounded by $\epsilon/3$. By definition of
$N^{(\infty)}_n$ in Section \ref{sec:model} and Little's formula,
$\expect N^{(\infty)}_n <\infty$ and, therefore,
\begin{equation}
\label{eq:ninftyn} \lim_{\psi\rightarrow \infty} \Pr \left
[N^{(\infty)}>\frac{\psi}{2}\right ]\rightarrow 0,
\end{equation}
uniformly for all $n>0$. Next, note that
$1[\theta_i^{0}>(1-\epsilon_1)n\expect [\tau_1-\tau_0]]\le
1[\theta_i^{0}>(1-\epsilon_1)\expect [\tau_1-\tau_0]]$ a.s., and
that for any fixed $m$,
\begin{equation*}
\Pr \left [\sum_{i=1}^{m}1[\theta_i^{0}>(1-\epsilon_1)n\expect
[\tau_1-\tau_0]]>\frac{\psi}{2}\right ]\le \Pr \left
[\sum_{i=1}^{m}1[\theta_i^{0}>(1-\epsilon_1)\expect
[\tau_1-\tau_0]]>\frac{\psi}{2}\right ]\downarrow
0\;\;\text{as}\;\;\text{$\psi\rightarrow \infty$,}
\end{equation*}
which by the monotone convergence theorem implies that the second
term in (\ref{eq:initstat1}) satisfies
\begin{equation}
\label{eq:sectermconv2} \lim_{\psi\rightarrow \infty}
\sum_{m=0}^{\infty}\Pr [N^{0}=m]\Pr \left
[\sum_{i=1}^{m}1[\theta_i^{0}>(1-\epsilon_1)n\expect
[\tau_1-\tau_0]]>\frac{\psi}{2}\right ]=0,
\end{equation}
uniformly for all $n>0$. Finally, by the Weak Law of Large Numbers,
for all $n$ large enough,
\begin{equation}
\label{eq:sectermpiece} \Pr [\tau_n-\tau_0\le (1-\epsilon_1)n\expect
[\tau_1-\tau_0]]\le \epsilon/9.
\end{equation}
Thus, the previous conclusion in conjunction with
(\ref{eq:sectermconv2}), (\ref{eq:ninftyn}) and (\ref{eq:initstat1})
implies that for an arbitrary $0<\epsilon <1$, there exist
$n_0<\infty$ and $\psi_0<\infty$ large enough such that for all
$n\ge n_0$
\begin{equation}
\label{eq:finsecterm2} \Pr \left
[N^{(\infty)}>\frac{\psi_0}{2}\right ]\le
\frac{\epsilon}{9}\;\;\text{and}\;\;\sum_{m=0}^{\infty}\Pr
[N^{0}=m]\Pr \left
[\sum_{i=1}^{m}1[\theta_i^{0}>(1-\epsilon_1)n\expect
[\tau_1-\tau_0]]>\frac{\psi_0}{2}\right ]\le \frac{\epsilon}{9}.
\end{equation}

Now, since $\expect B_i < \infty$, $\expect \theta_i <\infty$, there
exist $\beta_0$, $\delta_0$, such that
\begin{equation*}
\Pr [B_i >\beta_0] \le \frac{\epsilon}{3\psi_0}\;\;\text{and}\;\;\Pr
[\theta_i >\delta_0] \le \frac{\epsilon}{3\psi_0}.
\end{equation*}
Thus, the previous expressions in conjunction with
(\ref{eq:finsecterm2}), (\ref{eq:sectermpiece}) and
(\ref{eq:point2p1}) imply that for all large enough $n\ge n_0$
inequality (\ref{eq:point2}) holds for a chosen set
$S(\psi_0,\beta_0,\delta_0)$.

Next, we show that there exists $n_1 >0$ and $k>0$ such that for all
points $x\in S(\psi_0, \beta_0,\delta_0)$ and measurable sets
$A\subset S(\psi_0,\beta_0,\delta_0)$, the following inequality
holds
\begin{equation}
\label{eq:point1ineq1} P_{n_1} (x,A)\ge k R(A).
\end{equation}
Let $F_{\theta} (u)$ denote a cumulative distribution function of a
random duration $\theta$, i.e., $\Pr [\theta \le u]$. Furthermore,
select a small positive number $\eta$ such that for some chosen
$\Delta
> \delta_0$, $F_{\theta} (\Delta) -F_{\theta} (\delta_0) =\eta
>0$.
Next, for any $n_1$
\begin{equation}
\label{eq:point1main} P_{n_1} (x,A)\ge P_1 (x,\omega_0) P_{n_2}
(\omega_0,A),
\end{equation}
where $n_2 = n_1 - 1$. Let $x=(m,b_1,\dots, b_m,e_1,\dots ,e_m)\in
S(\psi_0, \beta_0,\delta_0)$. Then,
\begin{align}
P_1 (x,\omega_0) &\ge \Pr [\tau_1 -\tau_0 \ge \Delta , \text{all
$m$ requests depart in $(\tau_0,\tau_1)$}]\nonumber\\
&= \int_{\Delta}^{\infty} \prod_{i=1}^{m}\frac{F_{\theta} (u+e_i)
-F_{\theta} (e_i)}{1-F_{\theta} (e_i)} dF_a
(u),\label{eq:point1lowbd}
\end{align}
where $F_a (u)$ represents cumulative inter arrival distribution
of a renewal process $\{\tau_n\}$, i.e., $F_a (u)=\Pr
[\tau_1-\tau_0 \le u]$. Now, by applying lower bound
\begin{equation*}
\frac{F_{\theta} (u+e_i)-F_{\theta}(e_i)}{1-F_{\theta} (e_i)}\ge
F_{\theta} (\Delta)-F_{\theta} (\delta_0)=\eta
\end{equation*}
in (\ref{eq:point1lowbd}) we obtain
\begin{equation}
\label{eq:point1ineq2} P_1 (x,\omega_0)\ge \eta^{m}\Pr [\tau_1
-\tau_0
>\Delta]\ge \eta^{\psi_0} (1-F_a (\Delta)).
\end{equation}

Next, we derive a lower bound for $P_{n_2} (\omega_0,A)$ for some
$n_2$ large enough such that
\begin{equation}
\label{eq:wlln} \Pr [\tau_{n_2}-\tau_0 > \delta]\ge
1-\frac{\epsilon}{2}.
\end{equation}
Note that the condition imposed on $n_2$ in (\ref{eq:wlln}) is
possible due to the Weak Law of Large Numbers, since for any
$\epsilon >0$ and all $n_2$ large enough with $\delta_0 <
(1-\epsilon)\ex[\tau_{n_2}-\tau_0]$,
\begin{equation*}
\Pr [\tau_{n_2}-\tau_0>\delta_0]\ge \Pr
[\tau_{n_2}-\tau_0>(1-\epsilon)\expect [\tau_{n_2}-\tau_0]]\ge
1-\frac{\epsilon}{2},
\end{equation*}
Next, pick any $x'=(m',e'_1,\dots ,e'_{m'},b'_1,\dots,b'_{m'})\in
A$ where, without loss of generality, we assume that $e'_1\ge
e'_2\ge \dots \ge e'_{m'}$. Define $x'+dx'\triangleq
(m',e'_1+de'_1,\dots
,e'_{m'}+de'_{m'},b'_1+db'_1,\dots,b'_{m'}+db'_{m'})$ where
$de'_1, \cdots, de'_{m'}, db'_1,\cdots, db'_{m'}$ are
infinitesimal elements. Then, the transition probability into
state $(x', x'+d x')$ starting from $\omega_0$ can be bounded by
the probability of the event that there are exactly $m'$ arrivals
prior to $\tau_{n_2}$ whose arrival times are determined by $e'_1,
\cdots, e'_{m'}$, whose resource requirements are in
$(b'_1,b'_1+db'_1),\cdots, (b'_{m'},b'_{m'}+db'_{m'})$, and where
the rest of $n_2-m'$ arrivals are rejected since their
requirements exceed capacity $C$. Therefore,
\begin{align*}
&P_{n_2}(\omega_0,(x', x'+dx'))\ge \left \{ \prod_{j=1}^{m'}\Pr
[\theta_j
> e'_j]\Pr [B_j\in (b'_j,b'_j+db'_j)]\right \} \Pr
[B_1>C]^{n_2 - m'}\\
&\qquad\times\Pr \left [\bigcup_I \left \{\sum_{j=i_1}^{i_{2}-1}Y_j
\in (e'_1-e'_2,e'_1-e'_2+de'_1),\dots \sum_{j=i_{m'}}^{n_2-1}Y_j\in
(e'_{m'},e'_{m'}+de'_{m'})\right
\}\right ]\\
&\ge \Pr [B_1>C]^{n_2} \left \{ \prod_{j=1}^{m'}\Pr [\theta_j
> e'_j]\Pr [B_j\in (b'_j,b'_j+db'_j)]\right \}\\
&\qquad\times\Pr \left [\bigcup_I \left \{\sum_{j=i_1}^{i_2
-1}Y_j\in (e'_1-e'_2,e'_1-e'_2+de'_1),\dots
,\sum_{j=i_{m'}}^{n_2-1}Y_j\in (e'_{m'},e'_{m'}+de'_{m'})\right
\}\right ],
\end{align*}
where $I\triangleq\{0\le i_1<i_2<\dots <i_{m'}\le n_2 -1\}$ and
$Y_j$ are i.i.d. random variables equal in distribution to
inter-arrival times of the renewal process $\{\tau_n\}$, i.e.,
$Y_j\eqd \tau_{j+1}-\tau_j$. Now denote
\begin{align}
&r(m',e'_1,\dots ,e'_{m'}, b'_1,\cdots ,b'_{m'})\eqdef \left \{
\prod_{j=1}^{m'}\Pr [\theta_j
>
e'_j]\Pr [B_j\in (b'_j,b'_j+db'_j)]\right \} \nonumber \\
&\qquad \times\Pr \left [\bigcup_I \left \{\sum_{j=i_1}^{i_2
-1}Y_j \in (e'_1-e'_2,e'_1-e'_2+de'_1),\dots
,\sum_{j=i_{m'}}^{n_2-1}Y_j\in (e'_{m'},e'_{m'}+de'_{m'})\right
\}\right ],\label{eq:point1ineqdef}
\end{align}
and define probability distribution
\begin{equation}
\label{eq:defdistrpr} R(A)\eqdef {\mathcal V} \int_{x'\in A}
r(m',e'_1,\dots ,e'_{m'},b'_1,\dots ,b'_{m'}),
\end{equation}
where ${\mathcal V}$ is a normalization constant. Note that $R(A)$
is well-defined since
\begin{align*}
&\sum_{m'=0}^{\infty} \int_{\infty > e'_1 >\dots >e'_{m'} >0}
\int_{b'_1,\dots ,b'_{m'}\ge 0}^{\infty}r(m',e'_1,\dots
,e'_{m'},b'_1,\dots b'_{m'})\\
&\le \Pr [N^{(C)}_{0,n} \le \psi_0](\expect \theta_1)^{\psi_0} +\Pr
[N^{(C)}_{0,n} >\psi_0]\\
&\le \Pr [N^{(C)}_{0,n} \le \psi_0] (\expect \theta_1)^{\psi_0} +\Pr
[N^{(\infty)}_{n}
>\psi_0] < \infty.
\end{align*}
The previous inequalities, in conjunction with
(\ref{eq:point1ineqdef}), (\ref{eq:point1ineq2}) and
(\ref{eq:point1main}) imply that
\begin{equation}
\label{eq:point1} P_{n_2+1}(x,A)\ge \eta^{\psi_0} (1-F_a
(\Delta))\Pr [B_1>C]^{n_2} {\mathcal V}^{-1} R(A).
\end{equation}

Finally, it is left to show that there exists $K>0$ such that for
every initial distribution $P_0$, for all $n$ large and for any
measurable set $A\subset S(\psi_0,\beta_0,\delta_0)$
\begin{equation*}
P_n (A) \le K R(A)+\epsilon.
\end{equation*}
By (\ref{eq:wlln}), for all $n\ge n_2$
\begin{eqnarray*}
&P_n (A)& \le \Pr [X_{n} \in A, \tau_{n} -\tau_{0} >  \delta_0 ]+
\Pr
[\tau_{n} -\tau_{0} \le \delta_0]\\
&&\le \Pr [X_{n} \in A, \tau_{n} -\tau_{0} >  \delta_0 ]+\epsilon/2 \\
&&\le \int_{x'\in A}\left \{ \prod_{j=1}^{m'}\Pr [\theta_j
> e'_j]\Pr [B_j\in (b'_j,b'_j+db'_j)]\right \}\\
&&\qquad \times \Pr \left [\bigcup_I \left \{\sum_{j=i_1}^{i_2
-1}Y_j \in (e'_1 -e'_2,e'_1-e'_2+de'_1),\dots
,\sum_{j=i_{m'}}^{n_2-1}Y_j\in
(e'_{m'},e'_{m'}+de'_{m'})\right \} \right ]+\epsilon\\
&&=\int_{x'\in A}r(m',e'_1,\dots ,e'_{m'},b'_1,\dots
,b'_{m'})+\epsilon,
\end{eqnarray*}
where the second inequality follows from the fact that requests that
are active at $\tau_n$ must occur in the previous $\delta_0$ length
of time that are captured in $n_2$ renewal intervals
$[\tau_{n-n_2},\tau_{n-n_2+1}),\dots ,[\tau_{n-1},\tau_n)$ with
significant probability (greater than $1-\epsilon/2$). Thus, after
applying definition (\ref{eq:defdistrpr}), we obtain that for all
$n$ large
\begin{equation*}
P_n (A) \le {\mathcal V}^{-1} R(A)+\epsilon,
\end{equation*}
which, in conjunction with (\ref{eq:point1}) and (\ref{eq:point2}),
implies that the process $X_n$ satisfies conditions of the theorem
stated at the beginning of this section. Thus, there exists a unique
stationary distribution for the Markov chain $X_n$. Therefore, since
$Q_n^{(C)}$ defined in Section \ref{sec:model} is a functional of
the process $X_n$, it has a unique stationary distribution as well
implying the existence of the stationary blocking probability.

\qed

\section*{Acknowledgments}
The authors would like to thank Prof. Predrag Jelenkovi\'c for
valuable suggestions related to the possible generalizations of this
work.

\renewcommand{\baselinestretch}{1}
\small
\bibliography{anabib}

\begin{thebibliography}{10}

\bibitem{AHK94}
{\sc Asmussen, S., Henriksen, L.~F. and Kl\"{u}ppelberg, C.} (1994).
\newblock Large claims approximations for risk processes in a {M}arkovian
  environment.
\newblock {\em Stochastic Processes and their Applications\/} {\bf 54,} 29--43.

\bibitem{ASP05}
{\sc Asmussen, S. and Pihlg{\aa}rd, M.}
\newblock Loss rates for {L}{\'e}vy processes with two reflecting barriers.
\newblock {\em Mathematics of OR\/}.
\newblock to appear.

\bibitem{ATN72}
{\sc Athreya, K.~B. and Ney, P.~E.} (1972).
\newblock {\em Branching Processes}.
\newblock Springer-Verlag.

\bibitem{BON06}
{\sc Bonald, T.} (2006).
\newblock The {E}rlang model with non-{P}oisson call arrivals.
\newblock {\em Procedings of ACM/Sigmetric and Performance Conference,\/}
  276--286.

\bibitem{EKM97}
{\sc Embrechts, P., K.~C. and Mikosch, T.} (1997).
\newblock {\em Modelling Extremal Events for Insurance and Finance}.
\newblock Springer.

\bibitem{EMG80}
{\sc Embrechts, P. and Goldie, C.~M.} (1980).
\newblock On closure and factorization properties of subexponential and related
  distributions.
\newblock {\em J. Austral. Math. Soc.\/} {\bf Series(A) 29,} 243--256.

\bibitem{ERL17}
{\sc Erlang, A.~K.} (1917).
\newblock Solution of some problems in the theory of probabilities of
  significance in automatic telephone exchanges.
\newblock {\em Elektrotkeknikeren\/} {\bf 13,} 5--13.

\bibitem{GLM93}
{\sc Gazdzicki, P., Lambadaris, I. and Mazumdar, R.~R.} (1993).
\newblock Blocking probabilities for large multirate {E}rlang loss systems.
\newblock {\em Advances in Applied Probability\/} {\bf 25,} 997--1009.

\bibitem{GOK97}
{\sc Goldie, C.~M. and Kl\"uppelberg, C.} (1998).
\newblock Subexponential distributions.
\newblock In {\em A Practical Guide to Heavy Tails: Statistical Techniques for
  Analysing Heavy Tailed Distributions}. ed. M.~T. R.~Adler, R.~Feldman.
\newblock Birkh\"auser, Boston pp.~435--459.

\bibitem{HEL96}
{\sc Heyman, D.~P. and Lakshman, T.~V.} (1996).
\newblock Source models for {VBR} broadcast-video traffic.
\newblock {\em IEEE/ACM Transactions on Networking\/} {\bf 4,} 40--48.

\bibitem{JLS95a}
{\sc Jelenkovi\'c, P., Lazar, A. and Semret, N.} (1995).
\newblock Multiple time scales and subexponentiality in {MPEG} video streams.
\newblock {\em Technical report} CU/CTR/TR 430-95-36.
\newblock Columbia University.
\newblock http://www.ctr.columbia.edu/comet/publications.

\bibitem{JEL97c}
{\sc Jelenkovi\'{c}, P.~R.} (1999).
\newblock Subexponential loss rates in a {GI/GI/1} queue with applications.
\newblock {\em Queueing Systems, {S}pecial {I}ssue on {L}ong-{T}ailed
  {D}istributions\/} {\bf 33,} 91--123.

\bibitem{KAUF81}
{\sc Kaufman, J.~S.} (1981).
\newblock Blocking in a shared resources environment.
\newblock {\em \mbox{IEEE} Transactions on Communications\/} {\bf 29,}
  1474--1481.

\bibitem{KEL86}
{\sc Kelly, F.~P.} (1986).
\newblock Blocking probabilities in large circuit-switched networks.
\newblock {\em {A}dvances in {A}pplied {P}robability\/} {\bf 18,} 473--505.

\bibitem{KEL91s}
{\sc Kelly, F.~P.} (1991).
\newblock Loss networks.
\newblock {\em {A}nnals of {A}pplied {P}robability\/} {\bf 1,} 319--378.

\bibitem{KRM97}
{\sc Krishnan, K.~R. and Meempat, G.} (1997).
\newblock Long-range dependence in {VBR} video streams and atm traffic
  engineering.
\newblock {\em Performance Evaluation\/} {\bf 30,} 46--56.

\bibitem{LKT90}
{\sc Liu, L., Kashyap, B. R.~K. and Templeton, J. G.~C.} (1990).
\newblock On the ${GI}^{X}/{G}/\infty$ system.
\newblock {\em Journal of Applied Probability\/} {\bf 27,} 671--683.

\bibitem{LOU94}
{\sc Louth, G., Mitzenmacher, M. and Kelly, F.} (1994).
\newblock Computational complexity of loss networks.
\newblock {\em Theoretical Computer Science\/} {\bf 125,} 45--59.

\bibitem{LU06}
{\sc Lu, Y., Radovanovi\'c, A. and Squillante, M.} (2006).
\newblock Workforce management through stochstic network models.
\newblock {\em Proceedings of IEEE SOLI Conference\/}.

\bibitem{SEV57}
{\sc Sevastyanov, B.~A.} (1957).
\newblock An ergodic theorem for {M}arkov processes and its application to
  telephone systems with refusals.
\newblock {\em Theory of probability and its applications\/} {\bf 2,} 104--112.

\bibitem{TAK80}
{\sc Takacs, L.} (1980).
\newblock Queues with infinitely many servers.
\newblock {\em R.A.I.R.O. Recherche Operationnelle\/} {\bf 14,} 109--113.

\bibitem{WIT85}
{\sc Whitt, W.} (1985).
\newblock Blocking when service is required from several facilities
  simultaniously.
\newblock {\em AT$\&$T Technical Journal\/} {\bf 64,} 1807--1856.

\bibitem{WIL56}
{\sc Wilkinson, R.~I.} (1956).
\newblock Theory of toll traffic engineering in the {USA}.
\newblock {\em Bell System Technical Journal\/} {\bf 35,} 421--513.

\bibitem{WOL89}
{\sc Wolff, R.~W.} (1989).
\newblock {\em Stochastic Modeling and Theory of Queues}.
\newblock Prentice Hall.

\bibitem{ZAC91}
{\sc Zachary, S.} (1991).
\newblock On blocking in loss networks.
\newblock {\em {A}dvances in {A}pplied {P}robability\/} {\bf 23,} 355--372.

\end{thebibliography}
\end{document}